\documentclass[12pt]{article}
\usepackage{amsmath}
\usepackage{amssymb}
\usepackage{dsfont}
\usepackage{cite}
\newsymbol \blackbox 1004
\newcommand{\eh}{\hfill}\newlength{\sperr}

\newenvironment{proof}{{\settowidth{\sperr}{\bf\rm
Proof}%
\par\addvspace{0.3cm}\noindent\parbox[t]{1.3\sperr}
{\bf\rm P\eh r\eh o\eh o\eh f\eh }%
}}{\nopagebreak\mbox{}
$\blackbox$\par\addvspace{0.3cm}}

\def\nn{\nonumber}

\def\a{\alpha}

\def\g{\gamma}

\def\vk{\varkappa}

\def\Lam{\Lambda}
\def\s{\sigma}
\def\la{\lambda}
\def\om{\omega}

\def\vp{\varphi}

\def\ve{\varepsilon}
\def\wh{\widehat}
\def\wt{\widetilde}
\def\ov{\overline}

\def\p{\partial}

\def\BC{{\mathbb C}}

\def\BR{{\mathbb R}}
\def\BN{{\mathbb N}}

\def\cla{{\mathcal A}}
\def\clb{{\mathcal B}}
\def\clc{{\mathcal C}}

\def\cli{{\mathcal I}}

\def\cls{\mathcal{S}}

\def\cld{{\mathcal D}}

\newcommand{\E}{\mathrm{e}}
\newcommand{\I}{\mathrm{i}}

\newtheorem{Pa}{Paper}[section]
\newtheorem{Tm}[Pa]{{\bf Theorem}}

\newtheorem{Cy}[Pa]{{\bf Corollary}}
\newtheorem{Rk}[Pa]{{\bf Remark}}

\newtheorem{Ee}[Pa]{{\bf Example}}
\newtheorem{Dn}[Pa]{{\bf Definition}}

\newtheorem{Pn}[Pa]{{\bf Proposition}}

\title{On  new classes of explicit solutions   of Dirac, \\  dynamical Dirac and Dirac--Weyl  systems  \\ with non-vanishing at infinity  potentials,\\  their properties and applications}

\author{Alexander Sakhnovich}

\date{}
\begin{document}
\maketitle

\begin{abstract} Our GBDT version of B\"acklund-Darboux transformation is applied to the
construction of wide classes of new explicit solutions of self-adjoint and skew-self-adjoint Dirac systems,   dynamical Dirac and Dirac--Weyl  systems.
That is, we construct   explicit solutions of systems
with non-vanishing at infinity  potentials. In particular, the cases of steplike potentials and
power growth of potentials   are treated. It is essential (especially, for  dynamical case) that the generalised  matrix eigenvalues
are used in GBDT instead of the usual eigenvalues  (and those matrix eigenvalues are not necessarily
diagonal). The connection of Dirac--Weyl system with graphene theory is discussed.
Explicit expressions for Weyl--Titchmarsh functions are derived.

\end{abstract}


\vspace{0.2em}


\section{Introduction} \label{intro}
\setcounter{equation}{0}
In this paper we consider self-adjoint Dirac systems
\begin{align} &       \label{1.1}
y^{\prime}(x, z )=\I \big(z j+jV(x)\big)y(x,
z ), \quad
x \in \cli,
\\ &   \label{1.2}
j := \left[
\begin{array}{cc}
I_{p} & 0 \\ 0 & -I_{p}
\end{array}
\right], \hspace{1em} V(x)= \left[\begin{array}{cc}
0&v(x)\\  v(x)^* &0\end{array}\right],
 \end{align} 
and skew-self-adjoint Dirac systems
\begin{align} &       \label{1.3}
y^{\prime}(x, z )=\big(\I z j+jV(x)\big)y(x,
z ), \quad
x \in \cli,
 \end{align} 
where  $y^{\prime}:=\frac{\p}{\p x}y$, $\I$ stands for the  imaginary unit ($\I^2=-1$), $I_{p}$ is the $p \times p$ $(p\in \BN)$ identity
matrix,  $\BN$ stands for the set of positive integer numbers,  $v(x)$ is an $p \times p$ matrix function, and $\cli$ is some (finite or infinite)  interval
$(0\in \cli)$
on the real axis $\BR$.

We apply the obtained results to the study of the dynamical Dirac systems: 
$$\psi_x+\I j\big(\psi_\xi + V(x) \psi\big)=0 \quad {\mathrm{and}} \quad  \psi_x+\I j\big(\psi_\xi +\I  V(x)\psi\big)=0,$$
where $ \psi_x(x,\xi):=\frac{\p}{\p x}\psi(x,\xi)$ and $V$ does not depend on $\xi$. Then we turn to  Dirac--Weyl system (non-zero energy case) $$\psi_x=\I \s_3(-\psi_\xi +\I \om(x)\s_2\psi),$$ 
where $\s_k$ are Pauli matrices and $\om(x)=\ov{\om(x)}$.
Dirac--Weyl systems describe electron dynamics  and are  of interest in the study of graphene \cite{HaP, HoRoy, Mid}.

Systems \eqref{1.1}
and \eqref{1.3} play an essential role in  applications  to well known integrable wave equations and are also called {\it canonical, Dirac-type, Zakharov-Shabat or AKNS systems}.

Explicit solutions of systems \eqref{1.1} and \eqref{1.3} as well as explicit solutions of the corresponding direct and inverse
problems have been studied in detail for the case of $\cli=[0,\infty)$, potentials $V(x)$ tending to zero (as $x$ tends to infinity)
and rational Weyl functions and reflection coefficients (see, e.g., \cite{AlGo, FKRS-LAA, GKS2, GKS6, SaA2, ALS-JDE18, SaSaR, T0} and
references therein). The potentials corresponding to rational Weyl functions are called pseudo-exponential and the fact
that they tend to zero is stated in \cite[Proposition 4.1]{GKS6} for the self-adjoint case of system \eqref{1.1} and in \cite[Corollary 3.6]{FKRS-LAA}
for the skew-self-adjoint case of system \eqref{1.3}. In the present paper, we consider mostly both systems together (in the same statements).

The theory of explicit solutions in the case of potentials $V(x)$ which do not tend to zero  is of great interest.
See, for instance,  \cite{BKSh, HaP, LuMa1, LuMa2, SaA11} and references therein.
(See also interesting literature on the special cases of steplike potentials for 
other equations.)
However,  the case of potentials $V(x)$ which do not tend to zero is studied much less than the case
of potentials $V(x)$ tending to zero. 

We construct and study explicit solutions using our GBDT version of Darboux transformation.
B\"acklund-Darboux transformation is a well-known tool in the spectral theory and theory of  explicit solutions.
 (see, e.g.,  \cite{Ci, Gu, Mar, MS, Mi, ALS10, SaSaR, ZM} and numerous references therein). 
 See also related commutation methods \cite{D, Ge, GeT, T0}.
GBDT was first introduced in \cite{SaA2}, and 
a more general version of GBDT for first order systems rationally depending on the spectral parameter
was treated in \cite{ALS10, SaSaR} (see also some references therein).

The GBDT approach to Dirac systems is shortly discussed in section ``Preliminaries".
It is essential (for the dynamical case, in particular) that {\it generalised  matrix eigenvalues}
are used in GBDT instead of the usual eigenvalues  (and those matrix eigenvalues are not necessarily
diagonal).

Explicit solutions of Dirac systems with non-vanishing at infinity potentials
are considered in Section \ref{Tran}. Explicit solutions of dynamical Dirac and Dirac--Weyl systems
are constructed in Section \ref{DW}. The examples with steplike  potentials and potentials with
power growth are treated there as well. Finally, explicit expressions for Weyl--Titchmarsh functions are derived
in Section \ref{WTF}.

As usual, $\BR$ stands for the real axis and $\BC$ stands for the complex plane. 
The notation $\ov{a}$ means complex conjugate of $a$ and $A^*$ means complex conjugate transpose of the matrix $A$.
The open upper half-plane is denoted
by $\BC_+$ and the half-plane $\{z: \, \Im(z)>M>0\}$ is denoted by $\BC_M$.
By $\BC^{n\times k}$ ($\BR^{n\times k}$) we denote the class $n\times k$ matrices with complex-valued (real-valued) entries.
The inequality $S>0$ for some matrix $S$ means that $S$ is positive definite.
The notation $\Im(\a)$ stands for the imaginary part of matrix $\a$ (i.e.,  $\Im(\a)=\frac{1}{2\I}(\a-\a^*)$).
The space of square summable functions on $\cli$ is denoted by $L_2(\cli)$ and $L_2((0,\infty))=L_2(0,\infty)$.
By $L_2^{n\times k}(0,\infty)$ we denote the class of $n\times k$ matrix functions with the entries belonging
to $L_2(0,\infty)$.

\section{Preliminaries} \label{GBDT}
\setcounter{equation}{0}
GBDT of Dirac systems is determined by the initial system or, equivalently, by the potential $v(x)$ (or $V(x)$, see \eqref{1.2})
and by a triple of parameter matrices $\{A, S(0), \Pi(0)\}$ (see, e.g., \cite[Subsection 1.1.3]{SaSaR}).
Here, $A$ and $S(0)$ are $n \times n$ matrices $(n\in \BN)$, $\Pi(0)$ is an $n\times 2p$ matrix,
and the matrix identity
\begin{align}& \label{2.1}
AS(0)-S(0)A^*=\I \Pi(0) j^{\vk}\Pi(0)^*
\end{align}
holds, where 
\begin{equation} \label{2.2}
\vk=1 \quad {\mathrm{for \,\, system \,\, (\ref{1.1})}}; \quad \varkappa=0 \quad  {\mathrm{for  \,\, system \,\, (\ref{1.3})}}.
\end{equation}
{\it Our further formulas for the skew-self-adjoint case \eqref{1.3} somewhat differ from the self-adjoint case \eqref{1.1} as well.}

The so called Darboux matrix $w_A(x,\la)$ (which will be discussed later) is expressed via matrix functions $\Pi(x)$ and $S(x)$.
The matrix functions $\Pi(x)$  are determined by the  values $\Pi(0)$
and differential equations dual to Dirac systems:
\begin{align}& \label{2.3}
\Pi^{\prime}(x)=-\I A\Pi(x)j-\I\Pi(x)jV(x) \quad {\mathrm{for \,\, system \,\, (\ref{1.1})}}; 
\\ & \label{2.4}
\Pi^{\prime}(x)=-\I A\Pi(x)j-\Pi(x)jV(x)
\quad   {\mathrm{for  \,\, system \,\, (\ref{1.3})}}.
\end{align}
We partition $\Pi$ into $p\times p$ blocks $\Pi(x)=\begin{bmatrix}\Lam_1(x) & \Lam_2(x)\end{bmatrix}$
and rewrite \eqref{2.3} and \eqref{2.4} in the form
\begin{align}& \label{2.5}
\Lam_1^{\prime}=-\I A\Lam_1+\I\Lam_2 v^*, \quad \Lam_2^{\prime}=\I A\Lam_2-\I\Lam_1 v  \quad {\mathrm{for \,\, system \,\, (\ref{1.1})}}; 
\\ & \label{2.6}
\Lam_1^{\prime}=-\I A\Lam_1+\Lam_2 v^*, \quad \Lam_2^{\prime}=\I A\Lam_2-\Lam_1 v  \quad
   {\mathrm{for  \,\, system \,\, (\ref{1.3})}}.
\end{align}
The matrix functions $S(x)$  are determined by the  values $S(0)$
and matrix functions $\Pi(x)$ via equalities
\begin{align}& \label{2.7}
S(x)=S(0)+\int_0^x \Pi(r)j^{\vk+1}\Pi(r)^*dr,
\end{align}
where $\vk$ is given in \eqref{2.2}.  Relations \eqref{2.1}--\eqref{2.4} and \eqref{2.7} yield
\begin{align}& \label{2.1'}
AS(x)-S(x)A^*=\I \Pi(x) j^{\vk}\Pi(x)^* .
\end{align}
The transfer matrix function in Lev Sakhnovich form \cite{SaSaR, SaL1, SaL2} is given by the formula $w(z)=I_{2p}-\I  j^{\vk}\Pi^*S^{-1}(A-zI_n)^{-1}\Pi$,
and we consider the matrix function $w_A(x,z)$ introduced at each point $x$ of invertibility of $S(x)$ as the transfer matrix function:
\begin{align}& \label{2.8}
w_A(x,z)=I_{2p}-\I  j^{\vk}\Pi(x)^*S(x)^{-1}(A-zI_n)^{-1}\Pi(x),
\end{align}
where $\vk$ is introduced in \eqref{2.2}.
\begin{Pn} \label{Pn1} \cite{SaA2, SaSaR} Let \eqref{2.1} hold and let $w_A$ be given by \eqref{2.8}, where $\Pi(x)$ and $S(x)$ are given by \eqref{2.3} or
\eqref{2.4} and by \eqref{2.7}, respectively.

Then, for the self-adjoint system \eqref{1.1} we have
\begin{align}& \label{2.9}
\frac{d}{dx}w_A(x,z)=\I(zj+j\wt V(x))w_A(x,z)-\I w_A(x,z)\big(zj+j V(x)\big)  
\\ & \label{2.10}
\wt V(x)=V(x)+\I (\Pi(x)^*S(x)^{-1}\Pi(x) j-j\Pi(x)^*S(x)^{-1}\Pi(x)),
\\ &       \label{2.11}
\wt V(x)=\begin{bmatrix}0 & \wt v(x) \\ \wt v(x)^* & 0\end{bmatrix}, \quad
 \wt v(x)=v(x)-2 \I \Lam_1(x)^*S(x)^{-1}\Lam_2(x). \end{align}
For the skew-self-adjoint system \eqref{1.3} we have
\begin{align}& \label{2.12}
\frac{d}{dx}w_A(x,z)=(\I zj+j\wt V(x))w_A(x,z)- w_A(x,z)\big(\I zj+j V(x)\big)  
\\ & \label{2.13}
\wt V(x)=V(x)+\Pi(x)^*S(x)^{-1}\Pi(x) -j\Pi(x)^*S(x)^{-1}\Pi(x))j,
\\ &       \label{2.14}
\wt V(x)=\begin{bmatrix}0 & \wt v(x) \\ \wt v(x)^* & 0\end{bmatrix}, \quad
 \wt v(x)=v(x)+2 \Lam_1(x)^*S(x)^{-1}\Lam_2(x). \end{align}
\end{Pn}
\begin{Rk}\label{Rk1}
The matrix functions $w_A$ satisfying \eqref{2.9} or \eqref{2.12}
are called Darboux matrices for self-adjoint and skew-self-adjoint Dirac systems, respectively.
The corresponding Dirac systems with potentials $\wt V$ are called the transformed 
$($or GBDT-transformed$)$ Dirac 
systems determined by the triple $\{A,S(0), \Pi(0)\}$.
\end{Rk}
Proposition \ref{Pn1} is proved (see, e.g., \cite[Section 1.1.3]{SaSaR}) using the following relations.
\begin{Cy} \label{CyPiS} Let the conditions of Proposition \ref{Pn1} hold.\\
Then, in the self-adjoint case of system \eqref{1.1} and $\wt V$ given by \eqref{2.10} we have
\begin{align}& \label{2.15}
\big(\Pi^*S^{-1}\big)^{\prime}=\I j\Pi^*S^{-1}A+\I \wt V j \Pi^*S^{-1}.
\end{align}
In the skew-self-adjoint case of system \eqref{1.3} and $\wt V$ given by \eqref{2.13} we have
\begin{align}& \label{2.16}
\big(\Pi^*S^{-1}\big)^{\prime}=\I j\Pi^*S^{-1}A+j \wt V  \Pi^*S^{-1}.
\end{align}
\end{Cy}
The study of the invertibility of $S(x)$ is important, and we present below a simple but useful
statement  regarding this invertibility.
\begin{Pn}\label{Inv} Assume that $S(0)>0$. 

Then, $S(x)$ is invertible on $\cli=[0,\infty)$
in the self-adjoint case $($i.e., in the case of system \eqref{1.1}$)$, and $S(x)$ is invertible for all $x\in \cli= \BR$ in the skew-self-adjoint case
$($i.e., in the case of system \eqref{1.3}$)$.
\end{Pn}
\begin{proof}. In the self-adjoint case, we have $\vk=1$. Hence, \eqref{2.7} yields $S(x)>0$ for $x\geq 0$. The invertibility of $S(x)$
follows.

In the skew-self-adjoint case, we consider the matrix functions
\begin{align}& \label{2.17}
R(x):=\E^{-\I xA}S(x)\E^{\I xA^*}, \quad Q(x):= \E^{\I xA}S(x)\E^{-\I xA^*}.
\end{align}
(see \cite[Proposition 2.4]{SaADW}). Taking into account \eqref{2.7}, \eqref{2.1'} and \eqref{2.17}, we derive
\begin{align}\nn
R^{\prime}(x)&=\E^{-\I xA}\Pi(x)j\Pi(x)^*\E^{\I xA^*}-\I \E^{-\I xA}(AS(x)-S(x)A^*)\E^{\I xA^*}
\\ & \label{2.18}
=\E^{-\I xA}\Pi(x)(I_{2p}+j)\Pi(x)^*\E^{\I xA^*}\geq 0.
\end{align}
Since $R(0)=S(0)>0$, formula \eqref{2.18} implies that $R(x)>0$ for $x\geq 0$ (and so $S(x)$ is invertible
for $x\geq 0$. In the same way, we show that $Q^{\prime}(x)\leq 0$ and, hence, $Q(x)>0$ for $x\leq 0$.
Thus $S(x)$  is invertible
for $x\leq 0$ as well.
\end{proof}
\begin{Rk}\label{RkPos} The proof of Proposition \ref{Inv} shows  that $S(x)>0$ on $\cli$.
\end{Rk}
We will use the proof of Proposition \ref{Inv} and Remark \ref{RkPos} is some further considerations.
\section{Explicit solutions of the \\ transformed Dirac systems} \label{Tran}
\setcounter{equation}{0}
{\bf 1.} Instead of the trivial initial systems (i.e., systems with $v(x)\equiv 0$), 
which have been considered in previous research (see, e.g., \cite{AlGo, FKRS-LAA, GKS2, GKS6, SaSaR}),
let us consider initial systems
\eqref{1.1} and \eqref{1.3} with potentials
\begin{align}& \label{3.1}
v(x)=a\E^{2\I c x}I_p, \quad a\in \BC, \quad c\in \BR, \quad a\not=0.
\end{align}
\begin{Pn}\label{Lam} Let the  matrices $A$ and $Q$ $(A,Q\in \BC^{n\times n})$ satisfy relations
\begin{align}& \label{3.2}
AQ=QA, \quad Q^2=(A-cI_n)^2-|a|^2I_n \quad {\mathrm{for \,\, system \,\, (\ref{1.1})}}; 
\\ & \label{3.3}
AQ=QA, \quad Q^2=(A-cI_n)^2+|a|^2I_n 
\quad   {\mathrm{for  \,\, system \,\, (\ref{1.3})}}.
\end{align}
Let the matrices $f_1$ and $f_2$ belong $\BC^{n\times p}$ and
set
\begin{align}& \label{3.3'}
\Lam_1(0)=f_1+f_2, \quad \Lam_2(0)=f_3+f_4,
\end{align}
where $($in the case of the self-adjoint Dirac system$)$ $f_3$ and $f_4$ are given by
\begin{align}& \label{3.4}
f_3=(1/\ov{a})(Q+A-cI_n)f_1, \quad f_4=-(1/\ov{a})(Q-A+cI_n)f_2,
\end{align}
and in the case of the skew-self-adjoint Dirac system $f_3$ and $f_4$ are given by
\begin{align}& \label{3.5}
f_3=(\I /\ov{a})(Q+A-cI_n)f_1, \quad f_4=-(\I /\ov{a})(Q-A+cI_n)f_2.
\end{align}
Choose $S(0)=S(0)^*$ so that the identity \eqref{2.1} holds, where $$\Pi(0)=\begin{bmatrix}\Lam_1(0) & \Lam_2(0)\end{bmatrix}$$
is given by \eqref{3.3'}, \eqref{3.4} or \eqref{3.3'}, \eqref{3.5}.

Then, the triple $\{A, S(0), \Pi(0)\}$ determines GBDT of the corresponding Dirac system,
and the matrix functions $\Pi(x)=\begin{bmatrix}\Lam_1(x) & \Lam_2(x)\end{bmatrix}$ are expressed explicitly via formulas
\begin{align}& \label{3.6}
\Lam_1(x)=\E^{-\I c x}\big(\E^{\I  xQ}f_1+\E^{-\I  xQ}f_2\big), \quad \Lam_2(x)=\E^{\I c x}\big(\E^{\I  xQ}f_3+\E^{-\I  xQ}f_4\big),
\end{align}
where $\Lam_1$ and $\Lam_2$ satisfy \eqref{2.5} in the self-adjoint case and satisfy \eqref{2.6} in the skew-self-adjoint
case $($and $v$ is given by \eqref{3.1}$)$.
\end{Pn}
\begin{proof}. The proposition is proved by direct computation. We prove it for the self-adjoint case
(and the skew-self-adjoint case is proved in a similar way). According to \eqref{3.6} we have
\begin{align}& \label{3.7}
\Lam_1^{\prime}(x)=\I \E^{-\I c x}\big(\E^{\I  xQ}(Q-cI_n)f_1-\E^{-\I  xQ}(Q+cI_n) f_2\big), \\
& \label{3.8}
 \Lam_2^{\prime}(x)=\I \E^{\I c x}\big(\E^{\I  xQ}(Q+cI_n)f_3-\E^{-\I  xQ}(Q-cI_n)f_4\big).
\end{align}
From \eqref{3.1}, \eqref{3.6} and the first relation in \eqref{3.2} we have
\begin{align}& \label{3.9}
-\I A\Lam_1(x)+\I\Lam_2(x) v(x)^*=\I \E^{-\I c x}\big(\E^{\I  xQ}(-Af_1+\ov{a}f_3)+\E^{-\I  xQ}(-A f_2 +\ov{a}f_4)\big), \\
& \label{3.10}
\I A\Lam_2(x)-\I\Lam_1(x) v(x) =\I \E^{\I c x}\big(\E^{\I  xQ}(Af_3-af_1)+\E^{-\I  xQ}(Af_4-af_2)\big).
\end{align}
Comparing the right-hand sides of \eqref{3.7} and \eqref{3.9}, we see that they coincide
because (in view of \eqref{3.4}) the  terms containing $\E^{\I  xQ}$ and $\E^{-\I  xQ}$ coincide.
Thus, the first equality in \eqref{2.5} is satisfied. In order to show that the right-hand sides
of \eqref{3.8} and \eqref{3.10} coincide as well,  we use additionally the second equality
in \eqref{3.2}. Hence, the second equality in \eqref{2.5} is also valid.
\end{proof}
Clearly, the existence of matrices $Q$ satisfying \eqref{3.2} as well as matrices $Q$ satisfying \eqref{3.3} is essential,
and below we give an ``existence" proposition. Moreover, the proof of this proposition
presents a way to construct such $Q$.
\begin{Rk}\label{RkSqrt} The assertions including square roots $\sqrt{\ldots}$ further in the text are valid
for both values of the square roots $($if not stated otherwise$)$.
\end{Rk} 
\begin{Pn}\label{PnQ}  Let $a\not= 0$, $c\in \BR$ and an $n\times n$ matrix $A$ be given.
Then, 
if $\det\big(A-cI_n)^2-|a|^2I_n\big)\not=0$, there is an $n\times n$ matrix $Q$ such that \eqref{3.2} holds.
If $\det\big(A-cI_n)^2+|a|^2I_n\big)\not=0$, there is an $n\times n$ matrix $Q$ such that \eqref{3.3} holds.
\end{Pn}
\begin{proof}. We will construct $Q$ satisfying \eqref{3.2}, and $Q$ satisfying \eqref{3.3} is constructed in the same way.
Clearly, the statement of proposition is true for $n=1$. Consider the case, where $A$ is an $n\times n$ Jordan cell $(n\geq 2)$:
\begin{align}& \label{3.11}
A=\begin{bmatrix} \la & 1 & & 
\\  & \la & \ddots & \\
& & \ddots & 1 \\
& & & \la
\end{bmatrix}.
\end{align}
For this $A$, we construct an upper triangular Toeplitz matrix $Q$ satisfying the second equality in \eqref{3.2}.
Since $A$ is an  upper triangular Toeplitz matrix as well, it commutes with $Q$
(see, e.g., \cite{Com} on the properties of triangular Toeplitz matrices).
Thus, the first equality in \eqref{3.2} will be fulfilled for $Q$ automatically.

In order to construct $Q$, introduce the shift matrices
\begin{align}& \label{3.12}
 \cls_{i}:=\{\delta_{k-l+i}\}_{k,l=1}^n, \quad \cls_i \cls_j =\cls_{i+j}.
\end{align}
where $\delta_s$ is Kronecker delta, and $\cls_{i}=0$ for $i \geq n$.
Let us write down the representations
\begin{align}& \label{3.13}
A=\la I_n+\cls_1, \quad Q=q_0 I_n +q_1 \cls_1 + \ldots + q_{n-1}\cls_{n-1}.
\end{align}
According to \eqref{3.12} and \eqref{3.13}, we have
\begin{align} \label{3.14}
(A-cI_n)^2-|a|^2I_n=&\big((\la-c)^2-|a|^2) I_n+2(\la-c)\cls_1+\cls_2, \\
 \nn
Q^2=&q_0^2I_n+2q_0q_1\cls_1+(2q_0q_2+q_1^2)\cls_2
\\ & \label{3.15}
+ \sum_{i=3}^{n-1}\big(2q_0 q_i+q_1q_{i-1}+\ldots +q_{i-1}q_1\big)\cls_i.
\end{align}
Now, we set $q_0=\sqrt{(\la-c)^2-|a|^2}\not=0$  (recall Remark \ref{RkSqrt}) and choose consecutively the coefficients
$q_1, \, \ldots$ so that the coefficients before the shift matrices $\cls_i$ on the right-hand sides
of \eqref{3.14} and \eqref{3.15} coincide, where $q_i$ is uniquely determined by the
coefficients before $\cls_i$.

When $A$ is a Jordan matrix $J$,  we construct block diagonal  matrix $\wt Q$, each block of which is generated
by the corresponding Jordan block in a way described above. It is easy to see that
\eqref{3.2} holds for $A=J$ and $Q=\wt Q$. Finally, if $A=EJE^{-1}$, we set $Q=E\wt Q E^{-1}$
and equalities \eqref{3.2} for $A$ and $Q$ follow from \eqref{3.2} for $J$ and $\wt Q$.
\end{proof}
{\bf 2.}  Now, we can construct explicitly the potentials and fundamental solutions of the GBDT-transformed
Dirac systems.
\begin{Tm} \label{MTm} Let $a\not= 0$, $c\in \BR$ and an $n\times n$ matrix $A$ be given.\\
$($i$)$ Assume additionally that 
$\det\big(A-cI_n)^2-|a|^2I_n\big)\not=0$, choose $Q$ satisfying \eqref{3.2}, and let $w_A(x,z)$, $S(x)$ and $\Pi(x)$
be given explicitly by the formulas \eqref{2.8} $($with $\vk=1)$, \eqref{2.7} $($with $\vk=1)$ and \eqref{3.4}, \eqref{3.6}, respectively.

Then, a fundamental solution $\wt u(x,z)$ of the self-adjoint Dirac system
\begin{align} &       \label{1.1'}
y^{\prime}(x, z )=\I \big(z j+j\wt V(x)\big)y(x,
z ), 
 \end{align} 
where the potential $\wt V$ is given by \eqref{2.11} and \eqref{3.1}, has the form
\begin{align} &       \label{3.16}
\wt u(x,z)=w_A(x,z)u(x,z).
 \end{align} 
 Here, $u(x,z)$ is a fundamental solution of the initial system \eqref{1.1}
 and may be presented in the form
 \begin{align} &       \label{3.17}
u(x,z)=\E^{\I c x j}Z(z)\E^{\I x\zeta j}, \quad \zeta(z)=\sqrt{(z-c)^2-|a|^2},
\\ \label{3.18} &
Z(z)=\begin{bmatrix} aI_p & a I_p 
\\ (\zeta(z)-z+c)I_p & (-\zeta(z) -z +c)I_p
\end{bmatrix}.
 \end{align} 
 $($ii$)$  Assume  that 
$\det\big(A-cI_n)^2+|a|^2I_n\big)\not=0$, choose $Q$ satisfying \eqref{3.3}, and let $w_A(x,z)$, $S(x)$ and $\Pi(x)$
be given explicitly by the formulas \eqref{2.8} $($with $\vk=0)$, \eqref{2.7} $($with $\vk=0)$ and \eqref{3.5}, \eqref{3.6}, respectively.

Then, a fundamental solution $\wt u(x,z)$ of the skew-self-adjoint Dirac system
\begin{align} &       \label{1.3'}
y^{\prime}(x, z )= \big(\I z j+j\wt V(x)\big)y(x,
z ),  \end{align} 
where the potential $\wt V$ is given by \eqref{2.14} and \eqref{3.1}, has the form
\eqref{3.16}.
 Here, $u(x,z)$ is a fundamental solution of the initial system \eqref{1.3}
 and may be presented in the form
  \begin{align} &       \label{3.17'}
u(x,z)=\E^{\I c x j}Z(z)\E^{\I x\zeta j}, \quad \zeta(z)=\sqrt{(z-c)^2+|a|^2},
\\ \label{3.18'} &
Z(z)=\begin{bmatrix} \I aI_p & \I a I_p 
\\ (z-c-\zeta(z))I_p & (z -c+\zeta(z))I_p
\end{bmatrix}.
 \end{align} 
\end{Tm}
\begin{proof}.  Consider assertion (i). The existence of $Q$ follows Proposition \ref{PnQ}, the representations \eqref{3.4}--\eqref{3.6}
of the blocks $\Lam_1$ and $\Lam_2$ of $\Pi$ are shown in Proposition \ref{Lam}, and formula \eqref{3.16}
in immediate from Proposition \ref{Pn1}. It remains to prove that $u(x,z)$ of the form \eqref{3.17}, \eqref{3.18} is, indeed,
a fundamental solution of the initial system.

Let us partition $u$ into $p\times 2p$ blocks $u_1$ and $u_2$ and rewrite
Dirac system for a fundamental solution $u$ in the form
\begin{align} &       \label{3.19}
\frac{d}{dx}u_1(x,z)=\I z u_1(x,z)+\I v(x)u_2(x,z), \\
&     \label{3.20}
 \frac{d}{dx}u_2(x,z)=-\I z u_2(x,z)-\I v(x)^*u_1(x,z).
 \end{align} 
Relations \eqref{3.19}, \eqref{3.20} are dual to \eqref{2.5}, and $u$ may be constructed
using formula \eqref{3.6} for the blocks of $\Pi(x)$. However, it is simpler to
check that  \eqref{3.19} and \eqref{3.20} holds for $u$ of the form \eqref{3.17}
directly. Indeed, equalities \eqref{3.1} and \eqref{3.17}, \eqref{3.18} imply that
\begin{align} \nn       
\frac{d}{dx}u_1(x,z)=&\I z u_1(x,z)
\\ \nn &
+\I \E^{-\I c x}v(x)\begin{bmatrix}(-z+c+\zeta)\E^{\I x\zeta}I_p &(-z+c-\zeta)\E^{-\I x\zeta}I_p\end{bmatrix}
\\ \label{3.21} =&\I z u_1(x,z)+\I v(x)u_2(x,z).
 \end{align} 
The same equalities yield
\begin{align} \nn &
\frac{d}{dx}u_2(x,z)
\\ \nn
&=-\I z u_2(x,z)+\I \big(\E^{\I c x}v(x)^*/\overline{a}\big)
\\ & \nn \quad \times \begin{bmatrix}(z-c+\zeta)(-z+c+\zeta)\E^{\I x\zeta}I_p &(z-c-\zeta)(-z+c-\zeta)\E^{-\I x\zeta}I_p\end{bmatrix}
\\ \label{3.22} & =-\I z u_2(x,z)-\I v(x)^*u_1(x,z).
 \end{align} 
 Compare \eqref{3.19}, \eqref{3.20} with \eqref{3.21}, \eqref{3.22} to see that $u(x,z)$ given by \eqref{3.17}, \eqref{3.18} is the
 required fundamental solution.
 
 Assertion (ii) is proved quite similar to assertion (i).
\end{proof}
\begin{Rk} \label{RkS}  Let $\s(A)\cap \s(A^*)=\emptyset$ $($where $\s$ stands for spectrum$)$. Then,
$S(x)$ is uniquely recovered from the matrix identity \eqref{2.1'}. That is, \eqref{2.1'} may be used
instead of \eqref{2.7}.
\end{Rk}
\section{Dynamical Dirac and Dirac--Weyl systems} \label{DW}
\setcounter{equation}{0}
{\bf 1.} The motion of electron (in $\BR^2$)  in the presence of an electrostatic potential 
is often governed by the Dirac--Weyl system
\begin{align}& \label{I1}
\I \hbar v_F\left(\s_1 \frac{\p}{\p x}+\s_2 \frac{\p}{\p \xi}\right)\psi=\big(U(x,\xi)-E\big)\psi,
\end{align}
where $\hbar$ is the Planck constant, $v_F$ is the Fermi velocity, $U=\ov{U}$ and $E\in \BR$. 
Here, the matrices $\s_i$ are
Pauli matrices:
\begin{align}& \label{I2}
\s_1=\begin{bmatrix} 0 & 1\\  1 & 0 \end{bmatrix}, \quad
\s_2=\begin{bmatrix} 0 & - \I\\  \I & 0 \end{bmatrix}, \quad \s_3=\begin{bmatrix} 1 & 0 \\  0 & -1 \end{bmatrix}.
\end{align}
In the recent years, the interest in graphene has essentially stimulated the study of Dirac--Weyl system (see, e.g., \cite{HaP, HoRoy, Mid, SchUm, Stau}). In particular,
the important case when the scalar potential $U$ does not depend on the variable $\xi$ was recently studied,
for instance, in \cite{HaP, HoRoy, Mid} (see also some references therein). Assuming that $U$ does not depend on $\xi$ and multiplying both sides of \eqref{I1}
by $\frac{1}{\I \hbar v_F} \s_1$, we rewrite \eqref{I1} in an equivalent form
\begin{align}& \label{I3}
\psi_x=\I \s_3(-\psi_\xi +\I \om(x)\s_2\psi), \quad \psi_x:=\frac{\p}{\p x}\psi, \quad \om=\frac{E-U}{\hbar v_F}=\ov{\om}.
\end{align}
Some interesting cases of potentials generating explicit solutions
are studied  in \cite{HaP, HoRoy, Mid} using the separation of variables $\psi(x,\xi )=\E^{\I k\xi } \breve \psi(x)$, which transforms
\eqref{I3} into the system depending on one variable. Explicit solutions of the form $  \Psi(x)\E^{-\xi  A}g$, 
where $A$ are  $n\times n$ matrices, were constructed in \cite{SaADW} for the
case of so called pseudo-exponential potentials which tend to zero as $x$ tends to infinity. 

In this section, we will first consider  explicit solutions of dynamical Dirac systems
\begin{align}& \label{4.1}
\psi_x(x,\xi)+\I j\big(\psi_\xi(x,\xi) +\wt V(x) \psi(x,\xi)\big)=0,
\end{align}
and
\begin{align}
& \label{4.2}
\psi_x(x,\xi)+\I j\big(\psi_\xi(x,\xi) +\I \wt V(x)\psi(x,\xi)\big)=0
\end{align}
with non-vanishing at infinity potentials $\wt V$.
Then, we will consider Dirac--Weyl system \eqref{I3} as a particular case of system \eqref{4.2}. 

{\bf 2.} Our next theorem follows from Corollary \ref{CyPiS} and Theorem \ref{MTm}.
\begin{Tm} \label{DynTm} Let $a\not= 0$, $c\in \BR$ and an $n\times n$ matrix $A$ be given.

$($i$)$ Assume that 
$\det\big(A-cI_n)^2-|a|^2I_n\big)\not=0$, choose $Q$ satisfying \eqref{3.2}, and let $S(x)$ and $\Pi(x)$
be given explicitly by the formulas  \eqref{2.7} $($with $\vk=1)$ and \eqref{3.4}, \eqref{3.6}, respectively.

Then, the matrix function
\begin{align} &       \label{4.3}
\psi(x)=\Pi(x)^*S(x)^{-1}\E^{-\xi A}
\end{align} 
satisfies dynamical Dirac system \eqref{4.1}, where the potential $\wt V$ is given by \eqref{2.11} and \eqref{3.1}.

$($ii$)$  Assume  that 
$\det\big(A-cI_n)^2+|a|^2I_n\big)\not=0$, choose $Q$ satisfying \eqref{3.3}, and let $S(x)$ and $\Pi(x)$
be given explicitly by the formulas \eqref{2.7} $($with $\vk=0)$ and \eqref{3.5}, \eqref{3.6}, respectively.

Then, the matrix function $\psi$ of the form \eqref{4.3}
satisfies dynamical Dirac system \eqref{4.2}, where the potential $\wt V$ is given by \eqref{2.14} and \eqref{3.1}.
\end{Tm}
\begin{proof}. Differentiating $\psi$ of the form \eqref{4.3} and using  \eqref{2.15} we obtain
\begin{align}& \label{4.4}
\psi_x= \I j\big(\Pi(x)^*S(x)^{-1}A-\wt V(x)\Pi(x)^*S(x)^{-1}\big)\E^{-\xi A},
\end{align}
and assertion (i) follows. Differentiating $\psi$ of the form \eqref{4.3} and using  \eqref{2.16} we obtain
\begin{align}& \label{4.5}
\psi_x= \I j\big(\Pi(x)^*S(x)^{-1}A-\I\wt V(x)\Pi(x)^*S(x)^{-1}\big)\E^{-\xi A},
\end{align}
and assertion (ii) follows as well.
\end{proof}
{\bf 3.} Next, let $p=1$ and $\om(x)=-\I \wt v(x)=\ov{\om(x)}$. In this case, system \eqref{4.2} takes the
form of Dirac--Weyl system \eqref{I3}.
\begin{Cy} \label{CyDW} Let $p=1$, $c=0$, $a\not=0$,   
\begin{align}& \label{4.6}
\I a\in \BR, \quad \I A, \, \I Q, \, S(0)\in \BR^{n\times n}; \quad f_1, f_2 \in \BR^{n\times 1},
\end{align}
and assume that \eqref{3.3} in valid and \eqref{2.1} holds, where $\Pi(0)=\begin{bmatrix}\Lam_1(0) &
\Lam_2(0) \end{bmatrix}$ is given by \eqref{3.3'} and \eqref{3.5}. Let $\Pi(x)$ and $S(x)$
be given explicitly via \eqref{3.6} and \eqref{2.7} $($with $\vk=0)$, respectively.

Then, the matrix function $\psi(x)$ of the form \eqref{4.3} satisfies Dirac--Weyl system \eqref{I3},
where
\begin{align}& \label{4.7}
\om(x)=-\I \big(a+2\Lam_1(x)^*S(x)^{-1}\Lam_2(x)\big).
\end{align}
\end{Cy}
\begin{proof}. In view of Theorem \ref{DynTm} we need only to show that
\begin{align}& \label{4.8}
\om(x)=\ov{\om(x)}.
\end{align}
According to \eqref{3.5}, \eqref{3.6} and \eqref{4.6}, we have
\begin{align}& \label{4.9}
\Lam_1(x)\in \BR^{n\times 1}; \quad \I f_3, \, \I f_4, \, \I \Lam_2(x)\in \BR^{n\times 1}.
\end{align}
Moreover, relations \eqref{2.7}, \eqref{4.6} and \eqref{4.9} yield $S(x)\in \BR^{n \times n}$.
Hence, taking into account \eqref{4.7} and \eqref{4.9} we derive \eqref{4.8}.
\end{proof}
The simplest example corresponds to the case $n=1$.
\begin{Ee}\label{EeDW0}  Let $p=1$, $n=1$, $c=0$,
\begin{align}\label{4.25}&
a=\I r, \quad A=\I \la, \quad f_1=d, \quad f_2=1 \quad (r,\, \la,\,  d \,\in \BR), \\
\label{4.26} &
r\not=0,  \quad d\not=0, \quad \la > |r|, \quad Q=\I \mu, \quad  \mu=-\I \sqrt{r^2-\la^2}\in \BR
\end{align}
$($recall Remark \ref{RkSqrt}$)$. According to \eqref{3.6}, \eqref{4.25} and \eqref{4.26}, we have
\begin{align}\label{4.27}&
\Lam_1(x)=\E^{\mu x}+d\E^{-\mu x}, \quad \Lam_2(x)=\frac{\I}{r}\big((\mu-\la)\E^{\mu x}-d(\la+\mu)\E^{-\mu x}\big).
\end{align}
Relations \eqref{2.1'} and \eqref{4.27} yield
\begin{align}\label{4.28}&
S(x)=\frac{1}{2\la}\left(\Big(1+\frac{(\la-\mu)^2}{r^2}\Big)\E^{2\mu x}+4 d+d^2\Big(1+\frac{(\la+\mu)^2}{r^2}\Big)\E^{-2\mu x}\Big)\right).
\end{align}
Finally, formulas \eqref{4.7} and \eqref{4.25}--\eqref{4.27} imply that
\begin{align}\label{4.29}&
\om(x)=r+\frac{2}{r}\big((\mu-\la)\E^{2\mu x}-2d\la-d^2(\la+\mu)\E^{-2\mu x}\big)/S(x),
\end{align}
where $S(x)$ is given in \eqref{4.28}.
\end{Ee}

It is easy to see that $1+d^2+2d\geq 0$, $\frac{(\la-\mu)^2}{r^2}+d^2\frac{(\la+\mu)^2}{r^2}+2d\geq 0$, and the equalities
in the first and second inequalities  are achieved at different values of $d$. (Here, we take into account that $d(\la-\mu)(\la+\mu)=dr^2$
and $\mu\not=0$.) 
Hence, from
\eqref{4.28} we obtain $S(0)>0$. Thus, according to Remark \ref{RkPos}, the inequality $S(x)>0$ holds, and so $\om(x)$ is well defined.

It easily follows from \eqref{4.28} and \eqref{4.29} that $\om(x)$ is a {\it steplike potential}. Namely,
\begin{align}\label{4.30}&
\om(x)=r+\frac{4r{\la}(\mu-\la)}{r^2+(\mu-\la)^2}\big(1+O(\E^{-2\mu x})\big) 
\end{align}
in the cases $\mu>0, \,\, x\to \infty$ and $\mu<0, \,\, x\to - \infty$. We also have
\begin{align}\label{4.31}&
\om(x)=r-\frac{4r{\la}(\mu+\la)}{r^2+(\mu+\la)^2}\big(1+O(\E^{2\mu x})\big)
\end{align}
in the cases $\mu>0, \,\, x\to -\infty$ and $\mu<0, \,\, x\to  \infty$. 
\begin{Ee}\label{EeDW1}  Let us consider an example, where the matrices $A$ are not similar
to the diagonal matrices. Namely, we set $p=1$,  $c=0$, $n=2$,
\begin{align}& \label{4.10}
a=\I r, \quad A=\I\begin{bmatrix} \la & b \\ 0 & \la \end{bmatrix}, \quad f_1=\begin{bmatrix} d \\ 0  \end{bmatrix},
\quad f_2=\begin{bmatrix} 0 \\ 1  \end{bmatrix} \quad (r,\, \la,\, b,\, d \,\in \BR),
\\ &  \label{4.11}
r\not=0, \quad b\not=0, \quad d\not=0, \quad \la > |r|, \quad  \mu=-\I \sqrt{r^2-\la^2}\in \BR .
\end{align}
It is easily checked that
\begin{align}& \label{4.12}
Q=\I\begin{bmatrix} \mu & \frac{b \la}{\mu} \\ 0 & \mu \end{bmatrix}, \quad f_3=-\frac{\I d}{r}\begin{bmatrix} \la +\mu \\ 0  \end{bmatrix},
\quad f_4=\frac{\I}{r}\begin{bmatrix}  b \big((\la/\mu)-1\big) \\ \mu -\la  \end{bmatrix}.
\end{align}
Finally $($recall Remark \ref{RkS}$)$, the matrix function $S(x)$ is uniquely
defined by \eqref{2.1'}, which means that $S(x)=\{s_{ik}(x)\}_{i,k=1}^2$ is 
given by the equality
\begin{align}& \label{4.13}
2\la S(x)+\begin{bmatrix} b\big(s_{12}(x)+s_{21}(x)\big) & b s_{22}(x) \\ bs_{22}(x) & 0  \end{bmatrix}=\Lam_1(x) \Lam_1(x)^*(x)+\Lam_2(x)\Lam_2(x)^*.
\end{align}
In particular, relations \eqref{4.10}--\eqref{4.13} show that the conditions of Corollary~\ref{CyDW} are fulfilled.
\end{Ee}
{\bf 4.}  Let us study Example \ref{EeDW1} in greater detail.  From \eqref{4.12} we obtain
$$\E^{\pm \I x Q}=\E^{\mp  \mu x}\begin{bmatrix} 1 & \mp b\la x/\mu \\ 0 & 1  \end{bmatrix}.$$ 
Hence, relations \eqref{3.6}, \eqref{4.10} and \eqref{4.12} imply that
\begin{align}& \label{4.14}
\Lam_1(x)=\begin{bmatrix} \frac{b\la}{\mu}x\E^{\mu x}+d  \E^{-\mu x}\\ \E^{\mu x} \end{bmatrix}, \\
& \label{4.15}
\Lam_2(x)=\frac{\I}{r}\begin{bmatrix} 
(\mu-\la)\frac{b}{\mu}(\la x-1)\E^{\mu x}-d (\la +\mu) \E^{-\mu x}\\ (\mu - \la)\E^{\mu x} \end{bmatrix}.
\end{align}
Using equalities \eqref{4.13} and \eqref{4.14}, \eqref{4.15}, after some standard calculations
we obtain consecutively  the following relations:
\begin{align} \label{4.16}
s_{22}(x)&=\big({1}/{2\la}\big)\big(1+\big({(\mu-\la)}/{r}\big)^2\big)\E^{2\mu x}, \\
 \label{4.17}
s_{12}(x)&=s_{21}(x)=\frac{1}{2\la}\left(\frac{b\la}{\mu}\Big(1+\frac{(\mu-\la)^2}{r^2}\Big)x\E^{2\mu x}\right.
\\ \nn & \left.
\quad \qquad \qquad-
b\Big(\frac{(\mu-\la)^2}{r^2}\Big(\frac{1}{\mu}+\frac{1}{2\la}\Big)+\frac{1}{2 \la}\Big)\E^{2\mu x}+2d \right),
\end{align}
and
\begin{align}\nn
s_{11}(x)=&\frac{1}{2\la}
\left(\Big(\frac{b\la}{\mu}\Big)^2\Big(1+\frac{(\mu-\la)^2}{r^2}\Big)x^2\E^{2\mu x}
-\frac{b^2}{\mu}\Big(\Big(\frac{2\la}{\mu}+1\Big)\frac{(\mu-\la)^2}{r^2}+1\Big)\right.
\\  \label{4.18} & 
\times x\E^{2\mu x}
+b^2\Big(\frac{(\mu-\la)^2}{r^2}\Big(\frac{1}{\mu^2}+\frac{1}{\la\mu}+\frac{1}{2\la^2}\Big)+\frac{1}{2\la^2}\Big)\E^{2\mu x}
\\ & \nn \left.
+4bd\frac{\la}{\mu}x-2bd\Big(\frac{1}{\mu}+\frac{1}{\la}\Big)+d^2\Big(1+\frac{(\mu+\la)^2}{r^2}\Big)\E^{-2\mu x}\right).
\end{align}
Relations \eqref{4.16}--\eqref{4.18} provide explicit expressions for
\begin{align}& \label{4.19}
S(x)^{-1}=\frac{1}{\det S(x)}\begin{bmatrix} s_{22}(x) & -s_{12}(x) \\ -s_{21}(x) & s_{11}(x)  \end{bmatrix},
\end{align}
where we have
\begin{align} \nn
\det\big(S(x)\big)=&\frac{1}{4\la^2}\left(b^2\Big(\frac{1}{4\la^2}\Big(\frac{(\mu-\la)^4}{r^4}+1\Big)+\frac{(\mu-\la)^2}{r^2}\Big(\frac{1}{\mu^2}+\frac{1}{2\la^2}\Big)\Big)\E^{4\mu x}
\right.
\\ & \label{4.20} \left.
+\frac{2bd}{\mu}\Big(\frac{(\mu-\la)^2}{r^2}-1\Big)\E^{2\mu x}+4d^2\frac{\mu^2}{r^2}\right).
\end{align}
Thus, formulas  \eqref{4.14}--\eqref{4.20} give us expressions for the functions on the right-hand sides of \eqref{4.3} and \eqref{4.7}.
That is, they give us expression for the potential $\om$ and solution $\psi$ of the Dirac--Weyl system.

According to \eqref{4.16}, we have $s_{22}(0)>0$. It is easy to find sufficient conditions
for $\det S(0)>0$, which (together with $s_{22}(0)>0$) yields $S(0)>0$. In that case $S(x)$
is invertible (see Proposition \ref{Inv}). Therefore, the potential $\om$ and solution $\psi$ 
are well-defined on $\BR$.

The asymptotics of the potentials in Example \ref{EeDW1} is of interest. It  easily follows from \eqref{4.7} and \eqref{4.14}--\eqref{4.20}
and differs greatly
from the asymptotics of the pseudo-exponential potentials and steplike potentials. Recall that the value $\mu$ in   Example \ref{EeDW1} 
is given by \eqref{4.11}:
\begin{align} \label{4.21}
\mu=-\I \sqrt{r^2-\la^2}\in \BR \quad (\la >|r|).
\end{align}
\begin{Pn}. Let the conditions of Example \ref{EeDW1} hold. 

$(i)$ Assume that the branch of square root in \eqref{4.21} is chosen so that $\mu>0$.
Then, 
\begin{align} \label{4.22}
\om(x)=\frac{16 r(\la-\mu)\la^5\big(r^2+(\la-\mu)^2\big)}{\mu^2\big((\la-\mu)^4+r^4\big)+2(\la-\mu)^2r^2(2\la^2+\mu^2)}x^2\big(1+O(x^{-1})\big)
\end{align}
for $x\to \infty$, and 
\begin{align} \label{4.23}
\om(x)=r+O\big(\E^{2\mu x}\big) \quad {\mathrm{for}} \quad x \to -\infty .
\end{align}

$(ii)$ Assuming that the branch of square root in \eqref{4.21} is chosen so that $\mu <0$,
we have asymptotics \eqref{4.22} for $x \to -\infty$ and
\begin{align} \label{4.24}
\om(x)=r+O\big(\E^{2\mu x}\big) \quad {\mathrm{for}} \quad x \to \infty .
\end{align}
\end{Pn}
\section{Weyl-Titchmarsh functions} \label{WTF}
\setcounter{equation}{0}
Weyl--Titchmarsh theory of Dirac-type systems is an important and actively
developing domain (see, e.g., various general results and numerous references
in \cite{ClGe0, ClGe, FG-AS, SaA90, ALS-JDE18, SaSaR}). We note also that
references for the case of rational Weyl--Titchmarsh (Weyl) functions were adduced
in Introduction. Here, we consider explicit expressions for a wide new class of 
Weyl functions.

{\bf 1.} Let us consider Weyl functions of the self-adjoint Dirac systems \eqref{1.1'}
on the semi-axis $\cli=[0,\infty)$. We will use notations from \eqref{1.1'} and \eqref{1.3'} (in Theorem \ref{MTm})
instead of \eqref{1.1} and \eqref{1.3}, respectively, because later we turn from general-type potentials $\wt V$
to the class of potentials $\wt V$ considered in Theorem \ref{MTm}.

The notation $W(x,z)$ stands for the  normalised, by condition
\begin{align} \label{5.1}
W(0,z)=I_{2p},
\end{align}
fundamental solution of  the Dirac system  \eqref{1.1'}.
\begin{Dn} \label{DnWT}
A  $p\times p$ matrix function $\varphi$ such that
\begin{align} \label{5.2} &
\int_0^\infty \left[ \begin{array}{lr} I_p &  \I \varphi (z )^*
\end{array} \right]
  \Theta  W(x, z )^*
 W(x, z )\Theta ^*
 \left[ \begin{array}{c}
I_p \\ - \I \varphi (z ) \end{array} \right] dx < \infty , 
\\ \label{5.2'} &
 z 
\in {\BC}_+, \quad {\Theta}:=   \frac{1}{\sqrt{2}}       \left[
\begin{array}{cc} I_p &
-I_{p} \\ I_{p} & I_p
\end{array}
\right].
\end{align}
is called a Weyl function  of the Dirac  system \eqref{1.1'} on
$[0, \, \infty)$.
\end{Dn}
\begin{Rk} \label{RkWT}
There exists a unique Weyl function $\vp(z)$ of the system \eqref{1.1'} with locally summable potential $\wt V$. Moreover, this Weyl function is holomorphic
and $\Im\big(\vp(z)\big)>0$
 $($see, e.g., \cite[Subsection 2.1.1]{SaSaR} and the proof of \cite[Corollary 2.21]{SaSaR}$)$.
 \end{Rk}
 When the matrix function of the coefficients of linear fractional transformation admits a so called {\it realisation}, there is a standard
 way (see, e.g., \cite{GKS6}) to obtain a  realisation of the linear fractional transformation itself. We will need a corresponding proposition,
 which could be of independent interest.
 \begin{Pn}\label{Pnlft} Let $\clc_1,\clc_2\in \BC^{p\times n}, \quad \clb \in \BC^{n\times p}, \quad \cla \in \BC^{n\times n},$ and \\ $\cld \in \BC^{p\times p}$.
 Then, 
\begin{align} &\nn
 (\cld-\clc_2(\cla-z I_n)^{-1}\clb)(I_p-\clc_1(\cla-zI_n)^{-1}\clb)^{-1}
 \\ & \label{5.3} 
 =
 \cld+(\cld\clc_1-\clc_2)(\cla^{\times}-zI_n)^{-1}\clb, \quad \cla^{\times}:=\cla -\clb\clc_1.
 \end{align}
 \end{Pn} 
 \begin{proof}. It is well known that
 $$(I_p-\clc_1(\cla-zI_n)^{-1}\clb)^{-1}=I_p+\clc_1(\cla^{\times}-zI_n)^{-1}\clb.$$
 Hence, we have
 \begin{align} &\nn
 (\cld-\clc_2(\cla-z I_n)^{-1}\clb)(I_p-\clc_1(\cla-zI_n)^{-1}\clb)^{-1}
 \\ & \label{5.3+} 
 =\cld+\cld\clc_1(\cla^{\times}-zI_n)^{-1}\clb
 -\clc_2(\cla-zI_n)^{-1}\clb
\\ \nn & \quad 
 -
 \clc_2(\cla-zI_n)^{-1}\big((\cla-zI_n)-(\cla^{\times}-zI_n)\big)(\cla^{\times}-zI_n)^{-1}\clb,
 \end{align}
where $\cla^{\times}$ is given in \eqref{5.3}. The first equality in \eqref{5.3} follows from \eqref{5.3+}.
\end{proof} 
\begin{Tm} \label{TmWT} Let $a\not= 0$, $c\in \BR$ and an $n\times n$ matrix $A$ be given.
Assume that 
$\det\big(A-cI_n)^2-|a|^2I_n\big)\not=0$, choose $Q$ satisfying \eqref{3.2}, and let $S(x)$ and $\Pi(x)=\begin{bmatrix} \Lam_1(x) & \Lam_2(x)\end{bmatrix}$
be given explicitly by the formulas  \eqref{2.7} $($with $\vk=1)$ and \eqref{3.4}, \eqref{3.6}, respectively.
Assume additionally that $S(0)>0$.

Then, the Weyl function of the Dirac system \eqref{1.1'} on $[0,\infty)$, where $\wt V$ is given by
\eqref{2.11}, has the form
\begin{align} & \label{5.3'}
\vp(z)=\I\sqrt{2}\big(a-h(z)\big)^{-1}\big(\cld+(\cld\clc_1-\clc_2)(\cla^{\times}-zI_n)^{-1}\clb\big),
\\  \label{5.3!} &  \cld:=\frac{a+h(z)}{\sqrt{2}}I_p, \quad h(z):=\zeta(z)-z+c, \\
& \label{5.3!+}
 \clc_1:=\I\big(a-h(z)\big)^{-1}\big(\Lam_1(0)^*+\Lam_2(0)^*\big)S(0)^{-1}, \quad \cla^{\times}:=A -\clb\clc_1,
 \\ &  \label{5.3!!}
 \clb:=a\Lam_1(0)+h(z)\Lam_2(0), \quad \clc_2:=\frac{\I}{\sqrt{2}} \big(\Lam_1(0)^*-\Lam_2(0)^*\big)S(0)^{-1}.
\end{align}
$($The branch of $\zeta(z)$ in \eqref{5.3!} is chosen  so that $\zeta(z) \in \BC_+$ for $z\in \BC_+.)$
\end{Tm}
\begin{proof}.
According to Theorem \ref{MTm} and normalisation \eqref{5.1}, we have
\begin{align} \label{5.4}
W(x,z)=w_A(x,z)u(x,z)Z(z)^{-1}w_A(0,z)^{-1},
\end{align}
where $u$ and $Z$ are given by \eqref{3.17} and \eqref{3.18}, respectively. We choose the branch of $\zeta(z)$  in \eqref{3.17} and \eqref{3.18}
 so that $\zeta(z) \in \BC_+$ for $z\in \BC_+$. 

Consider the expression  \eqref{2.8} for $w_A(x,z)$. From (\ref{2.7})  (with $\vk=1$), it follows that
\begin{equation} \label{5.5}
\big(S(x)^{-1}\big)^{\prime}=-S(x)^{-1}\Pi(x) \Pi(x)^*
S(x)^{-1}.
\end{equation}
Hence, we obtain
\begin{equation} \label{5.6}
\int_0^{\infty}S(x)^{-1}\Pi(x) \Pi(x)^* S(x)^{-1}d x \leq
S(0)^{-1},
\end{equation}
that is, the entries of $\Pi^*S^{-1}$ belong to $L_2(0,\infty)$. 
In view of \eqref{3.6}, the entries of  $\E^{\I \zeta(z) x} \Pi(x)$
belong to $L_2(0,\infty)$ (when $\Im z>\|Q\|$) as well. Here, we used the equality
\begin{align} \label{5.7-}
\zeta(z)-(z-c)=-|a|^2\big(\zeta(z)+z-c\big)^{-1} ,
\end{align}
which yields $\Im(\zeta(z))>\Im(z)$ for $z\in \BC_+$.
The above-said  implies that
\begin{align} \label{5.7}
\E^{\I \zeta(z) x}w_A(x,z)\in L_2^{2p\times 2p}(0,\infty) \qquad (\Im(z)>\|Q\|).
\end{align}

Since the matrix $\Theta$ is unitary, relations \eqref{3.17}, \eqref{5.4} and \eqref{5.7} yield (in the case $\Im(z)>\|Q\|$)
\begin{align} \label{5.8}
W(x,z))\Theta ^*Y(z)\in L_2^{2p\times p}(0,\infty) \quad  {\mathrm{for}} \quad Y(z):=\Theta w_A(0,z)Z(z)\begin{bmatrix}I_p \\ 0 \end{bmatrix}.
\end{align}
Partition $Y(z)$ into the $p\times p$ blocks $Y_1(z)$ and $Y_2(z)$: $Y=\begin{bmatrix}Y_1 \\ Y_2\end{bmatrix}$, and compare \eqref{5.2} and \eqref{5.8} in order to see that
in the points of invertibility of $Y_1(z)$ (for $z$ satisfying $\Im(z)>\|Q\|$) the Weyl function $\vp$ has the form
\begin{align} \label{5.9}
\vp(z)=\I Y_2(z)Y_1(z)^{-1}.
\end{align}
In view of the equality \eqref{5.9}, of the definition of $Y$ in \eqref{5.8} and of Proposition \ref{Pnlft}, taking into account
the representations of $\Theta$, $w_A$ and $Z$ in \eqref{5.2'}, \eqref{2.8} and \eqref{3.18}, we derive \eqref{5.3'}--\eqref{5.3!!}
for $z$ satisfying $\Im(z)>\|Q\|$.
Since the Weyl functions are unique and holomorphic, we remove the condition  $\Im(z)>\|Q\|$ and the requirement
of the invertibility of $Y_1$.
\end{proof}
\begin{Rk} In the case of the pseudo-exponential potentials, we obtain not only explicit 
expressions for the Weyl functions but explicit expressions for the so called $\cla$-amplitudes
and for the inversion of the corresponding  convolution operators as well
$($see, e.g., \cite{AKAS}$)$. See also \cite{FG-AS} for the
$\cla$-amplitudes and corresponding structured operators for
Dirac-type systems. In the future, we plan to study the
$\cla$-amplitudes for the case of Dirac-type systems considered
in Theorem~\ref{TmWT}.
\end{Rk}
{\bf 2.} In a similar to Theorem \ref{TmWT} way, we construct explicitly  Weyl functions
of skew-self-adjoint Dirac systems. Recall the notation 
$$\BC_M=\{z: \, \Im(z)>M>0\}.$$
Similar to the self-adjoint case, the fundamental solution of the skew-self-adjoint Dirac system \eqref{1.3'} is denoted by $W(x,z)$ and
is normalised by \eqref{5.1}.
\begin{Dn} \label{DnW2} \cite[Section 3.1]{SaSaR}. A $p \times p$ matrix function $\vp$, which is  holomorphic in $\BC_M$
$($for some  $M>0)$ and satisfies the inequality
\begin{align}&      \label{5.11}
\int_0^{\infty}
\begin{bmatrix}
I_{p} & \vp(z)^*
\end{bmatrix}
W(x,z)^*W(x,z)
\begin{bmatrix}
I_{p} \\ \vp(z)
\end{bmatrix}dx< \infty , \quad z\in \BC_M,
\end{align}
is called a Weyl function of  the skew-self-adjoint Dirac system  \eqref{1.3'}.
\end{Dn}
Sometimes, it is more convenient to use the notion of the generalised Weyl function (GW-function). 
\begin{Dn}   \label{DnGW} \cite[Section 3.3]{SaSaR}. A   GW-function of  the system
\eqref{1.3'}, where $\wt V$ is locally bounded on $[0, \, \infty)$,
is a  $p \times p$ matrix function $\vp$ such that for some $M>0$ it is analytic in $\BC_M$ and the
inequalities
\begin{align}&      \label{5.12}
\sup_{x \leq \ell, \, z\in \BC_M}\left\| \E^{-\I z x}W(x,z)\begin{bmatrix}
I_{p} \\ \vp(z)
\end{bmatrix} \right\|<\infty
\end{align}
hold for each
$\ell < \infty$.
\end{Dn}
In particular, according to \cite[Proposition 3.28]{SaSaR} there is no more than one
GW-function for any system \eqref{1.3'} with  locally bounded potential $\wt V$.
\begin{Tm} Let $a\not= 0$, $c\in \BR$ and an $n\times n$ matrix $A$ be given.
Assume that 
$\det\big(A-cI_n)^2+|a|^2I_n\big)\not=0$, choose $Q$ satisfying \eqref{3.3}, and let $S(x)$ and $\Pi(x)=\begin{bmatrix} \Lam_1(x) & \Lam_2(x)\end{bmatrix}$
be given explicitly by the formulas  \eqref{2.7} $($with $\vk=0)$ and \eqref{3.5}, \eqref{3.6}, respectively.
Assume additionally that $S(0)>0$.

Then, the Weyl $($and GW-$)$ function of the Dirac system \eqref{1.3'} on $[0,\infty)$, where $\wt V$ is given by
\eqref{2.14}, has the form
\begin{align} & \label{5.16}
\vp(z)=\frac{1}{\I a}\big(\cld+(\cld\clc_1-\clc_2)(\cla^{\times}-zI_n)^{-1}\clb\big),
\\  \label{5.17} & \cld:=h(z)I_p, \quad  h(z):=z-c-\zeta(z),  \\
& \label{5.18}
 \clc_1:=\frac{1}{a_p}\Lam_1(0)^*S(0)^{-1}, \quad  \clc_2:=\I \Lam_2(0)^*S(0)^{-1},
 \\ &  \label{5.19}
\cla^{\times}:=A -\clb\clc_1, \quad \clb:=\I a\Lam_1(0)+h(z)\Lam_2(0). 
\end{align}
$($The branch of $\zeta(z)$ in \eqref{5.17} is chosen  so that $\zeta(z) \in \BC_+$ for $z\in \BC_+.)$
\end{Tm}
\begin{proof}.  According to Theorem \ref{MTm} and normalisation \eqref{5.1}, we have \eqref{5.4},
where $u$ and $Z$ are given by \eqref{3.17'} and \eqref{3.18'}, respectively. We choose the branch of $\zeta(z)$  in \eqref{3.17'} and \eqref{3.18'}
 so that $\zeta(z) \in \BC_+$ for $z\in \BC_+$. 

Consider again  the expression  \eqref{2.8} for $w_A(x,z)$, where $\vk=0$ in the present case.
Substitute the first equality in \eqref{2.17} into \eqref{2.8}:
\begin{equation} \label{5.13}
w_A(x,z)=I_{2p}-\I\Pi(x)^*\E^{\I x A^*}R(x)^{-1}(A-zI_n)^{-1}\E^{-\I x A}\Pi(x).
\end{equation}
In view of Remark \ref{RkPos}, we have $S(x)>0$. Hence, we derive from \eqref{2.17} that $R(x)^{-1}>0$ as well. 
The inequality \eqref{2.18} implies that $R(x)^{-1}$ is decreasing.

Instead of the equality \eqref{5.7-}, we have
\begin{align} \label{5.7!}
\zeta(z)-(z-c)=|a|^2\big(\zeta(z)+z-c\big)^{-1} ,
\end{align}
and so $\Im(\zeta(z))>\Im(z)-\ve$ for all $|z|$ such that $\Im(z)>1/\ve$ $(\ve>0)$.

Since $\Im(\zeta(z))>\Im(z)-\ve$ and $w_A$ admits representation \eqref{5.13}, where $R(x)^{-1}$ is decreasing and $R(x)^{-1}>0$,
relations \eqref{3.17'} and \eqref{5.4}  yield that for some $\g >0$,  $M>0$ and $\wt M>0$  we have
\begin{align} \label{5.8'}
\|W(x,z))Y(z)\|\leq \wt M \E^{-\g x} \quad (z\in \BC_M), \,\, {\mathrm{where}} \,\, Y(z):=w_A(0,z)Z(z)\begin{bmatrix}I_p \\ 0 \end{bmatrix}.
\end{align}
Moreover, for sufficiently large $M>0$ the matrix function $Y_1(z)$ is invertible and for some $\wh M>0$ we have
\begin{align} \label{5.14}
\| Y(z)^{-1}\|<\wh M \quad (z\in \BC_M).
\end{align}
Here, $Y_1(z)$ and $Y_2(z)$ are the $p \times p$ blocks of  $Y(z)$. Relations \eqref{5.8'} and \eqref{5.14} show that the GW-function and Weyl function
$\vp$ of the system \eqref{1.3'} is given by the formula:
\begin{align} \label{5.15}
\vp(z)=Y_2(z)Y_1(z)^{-1}.
\end{align}
The statement of the theorem follows from \eqref{5.15} and Proposition \ref{Pnlft}.
\end{proof}

{\bf Acknowledgments}  {This research    was supported by the
Austrian Science Fund (FWF) under Grant  No. P29177.}

\begin{flushright}
A.L. Sakhnovich,\\
Faculty of Mathematics,
University
of
Vienna, \\
Oskar-Morgenstern-Platz 1, A-1090 Vienna,
Austria, \\
e-mail: {\tt oleksandr.sakhnovych@univie.ac.at}

\end{flushright}

\end{document}